\documentclass{amsart}

\usepackage{amsmath}
\usepackage{amssymb}
\usepackage{url}
\usepackage{verbatim}

\newcommand{\nc}{\newcommand}
\nc{\gp}{\texttt{gp}}
\nc{\gap}{\texttt{gap}}
\nc{\Q}{{\bf Q}}
\nc{\Z}{{\bf Z}}
\nc{\Qp}{\Q_p}
\nc{\Zp}{\Z_p}
\nc{\QQ}{\Q_2}
\nc{\ZZ}{\Z_2}
\nc{\C}{{\bf C}}
\nc{\R}{{\bf R}}
\nc{\F}{{\bf F}}
\nc{\into}{\hookrightarrow}
\nc{\subeq}{\subseteq}

\nc{\Qbar}{\overline{\Q}}
\nc{\sbar}{\overline{s}}
\nc{\un}[1]{#1^{\mathrm{unr}}}
\nc{\g}{{gal}}
\nc{\emid}{\mid\kern -0.1ex\mid}

\nc{\OK}{{\mathcal{O}}_K}
\nc{\KK}{{\mathcal{K}_2}}
\nc{\G}{{\mathcal{G}_2}}
\nc{\tr}{\mathop{\rm Tr}}
\nc{\disc}{\mathop{\rm disc}}
\nc{\lcm}{\mathop{\rm lcm}}
\nc{\rd}{\mathop{\rm rd}}
\nc{\grd}{\mathop{\rm grd}}
\nc{\gms}{\textit{gms}}
\nc{\SC}{\mathrm{SC}}
\nc{\GSC}{\mathrm{GSC}}
\nc{\ead}{\mathrm{EAD}}
\nc{\syl}{\mathcal{S}_2}
\nc{\slope}{\textit{Slope}_{avg}}
\nc{\oord}{\mathrm{ord}}
\nc{\ord}{\oord_2}
\nc{\PSL}{PSL}
\nc{\SL}{SL}
\nc{\gal}{\mathop{\rm Gal}}
\nc{\Sym}{\mathop{\rm Sym}}
\nc{\Aut}{\mathop{\rm Aut}}
\nc{\cond}{\mathop{\rm cond}}
\nc{\Quot}{\mathop{\rm Quot}}
\nc{\cmmt}[1]{}
\nc{\gggg}{G(64,61)}
\nc{\gpsu}{PSU_3(3)}

\newtheorem{thm}{Theorem}[section]
\newtheorem{prop}[thm]{Proposition}
\newtheorem{cor}[thm]{Corollary}
\newtheorem{lemma}[thm]{Lemma}

\theoremstyle{remark}
\newtheorem{rem}[thm]{Remark}

\usepackage{array}
 \title{Number fields unramified away from $2$}
 \author{John W.\ Jones}

\begin{document}
  
\begin{abstract}
  Consider the set of number fields unramified away from $2$, i.e.,
  unramified outside $\{2, \infty\}$.  We show that there do not exist
  any such fields of degrees $9$ through $15$.  As a consequence, the
  following simple groups are ruled out for being the Galois group of
  an extension which is unramified away from $2$: Mathieu groups
  $M_{11}$ and $M_{12}$, $PSL(3,3)$, and alternating groups $A_j$ for
  $8<j<16$ (values $j\leq 8$ were previously known).
\end{abstract}

\maketitle

Let $\KK$ be the set of number fields $K\subset \C$ which are
unramified outside of the set $\{2, \infty\}$, i.e., fields with
discriminant $\pm 2^a$.  We say that such a field is unramified away
from $2$.  A field is in $\KK$ if and only
if its Galois closure is in $\KK$.  Accordingly, we let $\G$ be
the set of Galois groups fields in $\KK$ which are Galois over $\Q$.

Fields in $\KK$ and groups in $\G$ have been studied by several
authors \cite{tate, harb, brug, lesseni, marksau}.  In particular,
fields in $\KK$ of degree less than $9$ are fully understood, and a
variety of non-solvable groups have been shown to not lie in $\G$.
Here we extend these results for low degree fields.


The basic techniques used in the papers cited above are class field
theory, exhaustive computer searches of number fields with particular
discriminants, and discriminant bound arguments. 
In this paper, we will employ the third approach.  We use well-known lower
bounds for discriminants of number fields \cite{odlyzko-tables}.  Our upper
bounds for discriminants come from a study of higher ramification groups.
Preliminaries on discriminants of local fields are in
section~\ref{sec:local} with the main results in section~\ref{atmost15}.


\label{notation}


In general, we will use $K$ to denote a number field and $F$ to denote
a finite extension of $\Qp$, for some prime $p$.  Several notations
apply to both situations.  If $E$ is a finite degree $n$ extension of
$\Q$ or of $\Qp$, we let $(D_E)$ be the discriminant as an ideal over
the base, choosing $D_E$ to be a positive integer.  Then the root
discriminant for $E$ is $\rd(E) := D_E^{1/n}$.  We will denote the
Galois closure of $E$ over its base by $E^\g$.  Then, the Galois root
discriminant of $E$ is defined as $\grd(E) := \rd(E^\g)$.
%
We will use $\oord_p$ to denote the $p$-adic valuation on
$\Q$ such that $\oord_p(p)=1$.

When referring to Galois groups, we will use standard notations such
as $C_n$ for a cyclic group of order $n$, and $D_n$ for dihedral
groups of order $2n$.  Otherwise, we will use the $T$-numbering 
introduced in \cite{butler-mckay}, writing $nT_{j}$ for a degree $n$
field whose normal closure has Galois group $T_{j}$.


{\em Acknowledgments.} I would like to thank Dave Roberts for many
useful discussions. 

\section{Local fields} 
\label{sec:local}
Although we are particularly interested in extensions of $\QQ$, 
throughout this section, we work over $\Qp$ where $p$ is any prime.



\subsection{Slope}
If $F$ is a finite extension of $\Qp$,
then $D_{F}$ is a power of $p$. We define $c(F)$ to be the integer
such that $D_{F} = p^{c(F)}$.

Now assume $F$ is Galois over $\Qp$, with $G=\gal(F/\Qp)$.
We let $G^{\nu}$ be the higher ramification groups of $G$ in upper
numbering following the convention of \cite{jr-local-database}, which is
shifted by $1$ from \cite{serre-CL}.  In particular, $G=G^0$ and the
inertia subgroup is $G^1$.  We also define
$G^{\nu+} :=\bigcup_{\epsilon>0} \, G^{\nu+\epsilon}$.  Slopes of $F/\Qp$
are values $s$ where $G^{s+}\lneqq G^{s}$, i.e.
locations of jumps in the filtration.

Letting $\un{F}$ be the fixed field of $G^1=G^{0+}$ and
$F^{\mathrm{tame}}$ be the fixed field $G^{1+}$, we have that $\un{F}$
is the maximal unramified subextension of $F$ over $\Qp$ with
$\gal(\un{F}/\Qp) \cong G^0/G^{0+}$, and $F^{\mathrm{tame}}$ is the
maximal tame extension of $F$ over $\un{F}$ with
$\gal(F^{\mathrm{tame}}/\un{F}) \cong G^1/G^{1+}$.  Let
$f=[\un{F} : \Qp]=|G^0/G^{0+}|$
and $t=[F^{\mathrm{tame}} :  \un{F}]=|G^1/G^{1+}|$ be the
unramified and tame degrees respectively. These two integers
completely describe the only slopes $\leq 1$.

Slopes greater than $1$ correspond to wild ramification. The slope
content of $F/\Qp$ is then of the form $[s_1,\dots,s_m]_t^f$
where $f$ and $t$ are the unramified and tame degrees defined above,
and the $s_i$ are the wild slopes, sorted so that $s_i\leq
s_{i+1}$.  The ramification group
$G^{1+}$ is a $p$-group, and so for slopes $s>1$, the corresponding
quotients $G^{s}/G^{s+}$ are finite $p$-groups.  We repeat each $s_i$
with multiplicity $m_i$ where $p^{m_i}=|G^{s_i}/G^{s_i+}|$.
In particular, if $F/\Qp$ has slope content $[s_1,\dots,s_m]_t^f$,
then $|G|=p^m tf$.

Corresponding to the slope content $[s_1,\dots,s_m]_t^f$ is a
filtration on the Galois group which is just a slight modification of
the filtration discussed above.  For each wild slope $s>1$ with
multiplicity $k>1$, we refine the step $G^{s+}\leq G^s$ with index
$p^k$ into $k$ steps, each of degree $p$.  Taking fixed fields, we get
the tower
\begin{equation}\label{localtower}
  \Qp\subseteq \un{F}\subseteq F^{\mathrm{tame}}=F_0 \subseteq F_1
\subseteq F_2 \subseteq \dots \subseteq F
\end{equation}
where each extension $F_{i+1}/F_i$ has degree $p$.  For any finite
extension of local fields $F/E$, we define its average slope by
\begin{equation}  \label{slopedef}
  \slope(F/E) := \frac{c(F)-c(E)}{[F : \Qp] - [F :  E]}\, .
\end{equation}
For our tower~\eqref{localtower}, the average slopes give the wild
slopes in the slope content for $F/\Qp$, i.e.,
$s_i = \slope(F_i/F_{i-1})$ for $i\geq 1$
(see \cite{jr-local-database}).
So, the list of slopes can be discovered by working up through a
particular chain of subfields of $F/\Qp$.
%

From slope data $[s_1,\dots,s_m]_t^f$, we can compute the
discriminant, and hence the root discriminant of $F$.  Specifically,
$\rd(F) = p^{\gms_p(F)}$ where
\begin{equation}\label{gmsdef} 
\gms_p(F) = \frac{c(F)}{[F : \Qp]} = 
\frac{1}{p^m} \left(\sum_{i=1}^m p^{i-1} (p-1) s_{i} + 
 \frac{t-1}{t}\right)\, .
\end{equation}
Note, the unramified degree $f$ does not enter into
formula~\eqref{gmsdef}.  For the remainder
of this paper, we will omit $f$ from the slope content of an extension
and write simply $[s_1,\dots,s_m]_t$.  Since formula~\eqref{gmsdef} is
already a function of the slope content, we will also use it to define
$\gms_p([s_1,\dots,s_m]_t)$.  When comparing possible slope contents
$\alpha$ and $\beta$, we say that $\alpha$ is an upper bound for
$\beta$ if $\gms_p(\alpha)\geq \gms_p(\beta)$.
\begin{rem}
  The notation $\gms$ stands for {\em Galois mean slope}, so named
  because it is a weighted average of slopes for a Galois extension.
  The terminology is similar to our use of {\em average slope},
  denoted by $\slope$, which is also a weighted sum of slopes from a
  Galois extension.  However, we will not be making use of this latter fact
  here.
\end{rem}

\subsection{Composita}
\label{sec:composita}
If we start with a global field $K$, we can compute $\grd(K)$
locally.  We decompose $K\otimes \Qp\cong
\prod_{i=1}^g K_{p,i}$ as a product of finite extensions of $\Qp$.
The algebra $K^\g\otimes\Qp$ is a product of copies of $K^\g_p:=
(K_{p,1})^\g \cdots (K_{p,g})^\g$, the compositum of the Galois
closures of the $K_{p,i}$.
Picking a prime for $K^\g$ above each prime $p$, we let $\gms_p(K) :=
\gms_p(K_p^\g)$, and then
\[ \grd(K)=\rd(K^\g)=\prod_p p^{\gms_p(K)}\, .\]
Naturally, in this product, the factor for each unramified prime is
$p^0=1$.

An important, and somewhat subtle problem, then is to determine
$\gms_p$ for the compositum of fields $K_{p,i}^\g$.
Proposition~\ref{composprop} below gives reasonable bounds on $\gms_p$
for a compositum.
Given a slope content $\alpha=[s_1,\dots,s_m]_t$ and a rational number
$s>1$, we write $m_s(\alpha)$ for the multiplicity of $s$ in $\alpha$,
i.e., the number of $s_i$ equal to $s$.  Similarly, we write $m_{\geq
  s}(\alpha)$ for the number of slopes $s_i\geq s$.
The following proposition is a straightforward consequence of
Herbrand's theorem \cite{serre-CL}.
\begin{prop}\label{composprop}
  Suppose $F_1$ and $F_2$ are finite Galois extensions of $\Qp$, 
  with slope contents $\alpha_1$
  and $\alpha_2$.
  Let 
$\beta$
be the slope
  content of the compositum $F_1F_2$.  Then,
  \begin{enumerate}
  \item for all $s>1$, $m_s(\beta)\ge \max(m_s(\alpha_1), m_s(\alpha_2))$;
  \item for all $s>1$, $m_{\geq s}(\beta) \leq
    m_{\geq s}(\alpha_1)+m_{\geq s}(\alpha_2)$.
  \end{enumerate}
  Moreover, the tame degree for $F_1F_2$ is the least common multiple
  of the tame degrees of $F_1$ and $F_2$.
\end{prop}

Given two finite Galois extensions $F_1$ and $F_2$ of $\Qp$,
Proposition~\ref{composprop} gives upper bounds for the slope
content of $F_1F_2$, and hence for $\gms_p(F_1F_2)$, which are easy to
compute.  Namely, one
combines the tame degrees as described in the proposition, and just
concatenates (and sorts) the lists of wild slopes.  To bound the
slope content of the compositum subject to an upper bound on the number
of wild slopes, one removes slopes from the combination which occur in
the slope contents of both fields, starting with the smallest such
slopes.

For example, given slope contents $[2,3,7/2]_9$ and $[2,3,4]_{15}$, 
an upper bound for the slope content of the compositum is
$[2,2,3,3,7/2,4]_{45}$.  The maximal combinations with $5$ and $4$ wild
slopes are $[2,3,3,7/2,4]_{45}$ and $[2,3,7/2,4]_{45}$ respectively.
One cannot have a combination with less than $4$ wild slopes in this
case by Proposition~\ref{composprop}.  We will refer to this process
as computing the {\em crude upper bound} for slope content.  In
some cases, one can certainly obtain better bounds by using more
knowledge of the fields involved, e.g., Proposition~\ref{allquartic}
below.

\subsection{Individual slope bounds}

Our first lemma follows from basic facts about ramification
\cite[Chap. 3]{serre-CL}, and some simple algebra to translate
statements from discriminant exponents to slopes.
\begin{lemma} \label{appslope}
  If $F \supseteq E \supseteq \Qp$ are finite extensions where $F/E$
  is totally ramified of degree $p^n$, and $[E :  \Qp]=ef$
  where $f$ is the residue field degree for $E/\Qp$, then 
  \[ c(F) = p^n \cdot c(E) + f\nu \]
  where $\nu$ is an integer, $ep^n \leq \nu \leq p^n-1+ n ep^n$.
  Moreover, the average slope for $F/E$ equals
  \begin{equation} \label{appslFE} 
 \slope(F/E) = \frac{c(E)}{[E : \Qp]} + \frac{\nu}{(p^n-1)e}\, .
\end{equation}
\end{lemma}
\begin{rem}
  In Lemma~\ref{appslope}, given a field $E$, there
  will exist an extensions $F$ of degree $p^n$ satisfying both
  extremes of the inequalities for $\nu$.  If $\pi$ is a uniformizer
  for $E$, one can use $x^{p^n}+\pi x + \pi$ and $x^{p^n} + \pi$ to
  define extensions achieving the lower and upper bounds respectively.
\end{rem}

We now apply Lemma~\ref{appslope} to bound the average slopes in a tower.

\begin{lemma} \label{upperplem}
  Given a tower of finite extensions
  \[
  \Qp\subseteq F_0 \subseteq F_1 \subseteq F_2 \subseteq \dots
  \subseteq F_m= F
  \]
  where $[F_0:\Qp]=ef$, with $f$ being the residue field degree,
  $p\nmid e$,  and  each $F_i/F_{i-1}$ totally ramified of degree $p$,
  then for all $i$,   $1\leq i \leq m$,
  \[ \slope(F_i/F_{i-1}) \le i+\frac{p}{p-1}\, .\]
\end{lemma}
\begin{proof}
We will abbreviate $\slope(F_i/F_{i-1})$ by $S_i$.
Applying Lemma~\ref{appslope}, we get the recursions
\begin{equation} \label{crecursion}
  c(F_i) = p c(F_{i-1}) + f\nu_i
\end{equation}
and
\begin{equation}\label{srecursion}
  S_i = \frac{c(F_{i-1})}{p^{i-1}ef} +
  \frac{\nu_i}{(p-1)p^{i-1}e}   \,.
\end{equation}
Here $\nu_i$ is the value of $\nu$ in Lemma~\ref{appslope} for the
extension $F_i/F_{i-1}$.
Note that in tower, we have
\begin{equation}\label{dif1}
  S_{i+1}-S_i =
  \frac{c_{i+1}-c_i}{(p-1)p^{i}ef}-\frac{c_i-c_{i-1}}{(p-1)p^{i-1}ef}
  = \frac{\nu_{i+1}-\nu_i}{(p-1)p^{i}e}
\end{equation}
From equations~\eqref{crecursion} and \eqref{srecursion}, it is clear
that the sequences of 
discriminant exponents, $c_i$, and average slopes, $S_i$, are each
bounded by the 
corresponding sequences where we use the upper bound for each
$\nu_i \leq p-1+p^ie$.  For the sequence of $S_i$ where $\nu_i$ is
maximal for all $i$,
\begin{equation} \label{maxram}
  S_{i+1}-S_i 
  \leq \frac{(p-1+p^{i+1}e)-(p-1+p^ie)}{(p-1)p^{i}e} = 1
\end{equation}
So $S_i \leq S_1+i-1$, and it is easy to check from
equation~\eqref{srecursion} and the bound for $\nu_1$
that $S_1 \leq 1+p/(p-1)$, giving the result.
\end{proof}

\begin{rem}
  It is not always the case that in a sequence of slopes,  $s_i
  \leq s_{i-1}+1$ for all $i$.  For example, there are extensions of
  $\QQ$ with Galois group $8T_{23}$ with slope content $[4/3, 4/3,
  3]_3$ (see \cite{jr-2adicoctics}).
\end{rem}

We have two applications of Lemma~\ref{upperplem}.  First, we can
apply it directly to the tower in display~\eqref{localtower} to get bounds on
the wild slopes of a Galois extension $F/\Qp$.
\begin{prop}\label{upperp}
  Let $F/\Qp$ be a Galois extension with slope content
  $[s_1,\dots,s_m]_t$.  Then for $1\leq i \leq m$, 
\[ s_i \leq i+\frac{p}{p-1}\, . \]
\end{prop}

Specializing Proposition~\ref{upperp} to the case $p=2$, we have the
following.
\begin{cor} \label{upperslopes}
  If a Galois extension $F/\QQ$ has slope content $[s_1,\dots,
  s_m]_t$, then $s_i\leq i+2$ for $1\leq i\leq m$. 
\end{cor}
\begin{rem}
  The bound in Corollary~\ref{upperslopes} is achieved by a cyclic
  extension of degree $2^k$ over $\QQ$ given by $\QQ(\zeta_{2^k} +
  \zeta_{2^k}^{-1})$.
\end{rem}

We now apply Lemma~\ref{upperplem} to a non-Galois extension, and
slopes of its Galois closure.
\begin{prop}\label{slopebound12}
  If $F$ is a finite extension of $\Qp$, then all slopes $s$ for $F^\g/\Qp$
  satisfy 
  \[s \leq \frac{p}{p-1}+\oord_p([F:\Qp])\, .\]
\end{prop}
\begin{proof}
  Let $\sbar$ be the largest slope for $F^\g/\Qp$, so we need to show
  that $\sbar \leq \frac{p}{p-1}+\oord_p([F:\Qp])$.  This is clear if
  $\sbar\leq 1$, so we can assume that
  $\sbar>1$, i.e. there is wild ramification.

  Let $G =\gal(F^g/\Qp)$ and let $H$ be the subgroup fixing $F$.
  From \cite[\S3.6]{jr-local-database}, the extension $F/\Qp$ has a
  distinguished chain of subfields corresponding to subgroups $HG^s$;
  we will denote the fixed field of $HG^s$ by $F_s$, and define
  $F_{s+}$ analogously.
  For values of $s$ where $HG^s\neq HG^{s+}$, $s=\slope(F_{s+}/F_{s})$.
  Since $G^{\sbar+}$ is trivial and $H$ cannot contain a non-trivial
  normal subgroup of $G$, $HG^{\sbar}\neq H=HG^{\sbar+}$.  Hence,
  $\sbar=\slope(F/F_{\sbar})$.

  Among extensions $F/\Qp$ of a given degree, it is clear
  geometrically from \cite[\S3.6]{jr-local-database}, or algebraically
  from Lemma~\ref{appslope}, that the value of $\sbar$ for any
  extension is bounded by it value for an extension having
  intermediate fields of index $p^j$ for all $0\leq j\leq
  \oord_p([F:\Qp])$.  So, we can apply
  Lemma~\ref{upperplem} to obtain
  $\sbar \leq \frac{p}{p-1}+\oord_p([F:\Qp])$.
\end{proof}
\begin{rem}
  Proposition~\ref{slopebound12} will be applied below to extensions
  of $\QQ$ of degrees $12$ and $14$, showing that the Galois closure
  in each case has wild slopes bounded by $4$ and $3$ respectively.
\end{rem}
\begin{rem}
  The proof of Proposition~\ref{slopebound12} shows that the extension
  $F/\Qp$ must have certain intermediate fields, including a subfield
  corresponding the the largest slope for $F^g/\Qp$.  A nice
  illustration of this comes from sextic extensions of $\QQ$.  If the
  extension is wildly ramified, then $F/\QQ$ must have a cubic
  subfield.  Checking the appropriate table at
  \cite{web-local-database}, we see that there is exactly one sextic
  extension of $\QQ$ which does not have a cubic subfield,  but it is
  not wildly ramified.  For fields which are wildly ramified, the
  slope of $F/F_3$ where $F_3$ is the cubic subfield is the largest
  slope for $F^\g$.  For the field $F$ with no cubic subfield, there
  is tame ramification and $F$ has a quadratic subfield $F'$ so that
  $F/F'$ corresponds to the  maximum slope of $1$ for $F^\g$.
\end{rem}

\section{Number fields of degree less than $16$}
\label{atmost15}
In sections~\ref{deg9-11} and \ref{deg12}, we prove the following
theorem.
\begin{thm}\label{mainthm}
  There do not exist any degree $n$ extensions of $\Q$ which are
  unramified away from $2$ where $9\leq n \leq 15$.
\end{thm}
We consider each degree $n$, and within each degree we
consider the possible Galois groups among the transitive subgroups of
$S_n$.  To minimize the
number of cases we need to consider in detail, we note that if $G$ is
the Galois group of $K^\g\in \KK$ where $[K :  \Q]>8$, then $G$ must
satisfy the following two properties:
\begin{enumerate}
  \item $|G|$ is a multiple of $2^4$;
\item all proper quotients of $G$ are in $\G$.
\end{enumerate}
The first property is a consequence of Theorem~2.23 of \cite{harb};
the second is clear.
Progressing successively through degrees, there will only be a small
number of groups which satisfy both conditions.  For
reference, we state here previously known results of groups which are
not in $\G$ based on \cite{tate, harb, marksau, brug, lesseni}.  They
provide the starting point for applying property~(2) above.
Suppose $K\in \KK$ and $G\in\G$.  Then,
\begin{enumerate}
\item $[K : \Q]\neq 3$, $5$, $6$, $7$;
\item if $[K : \Q]\leq 8$, then $\gal(K^\g/\Q)$ is a $2$-group;
\item if $|G|<272$, then $G$ is a $2$-group;
\item $G\neq \PSL_2(2^j)$ for $j\ge 1$;
\item if $G$ is a $2$-group, then $G$ can be generated by two
  elements, one of which is $2$-torsion.
\end{enumerate}
Mark{\v{s}}a{\u\i}tis's result \cite{marksau} carries even more
information.  If $G_{\Q,2}$ is the Galois group of the maximal
extension of $\Q$ unramified away from $2$, he shows that the maximal
pro-$2$ quotient of $G_{\Q,2}$ is the pro-$2$ completion of the free
product $\Z *C_2$.

For lower bounds on root discriminants, we will refer to 
Table~\ref{uncond}.
\begin{table}[htb]
\caption{Unconditional root discriminant bounds.\label{uncond}  A
  field $K$ with $[K : \Q]\ge n$ has $\rd(K)$ greater than or equal to
  the given value.  If $K\in\KK$, then $\gms_2(K)$ is greater than or
  equal to the given bounds.  }
\[ \begin{array}{|c|r|r|c|c|r|r|}
\cline{1-3} \cline{5-7}
\gms_2\text{ for }\KK&
\multicolumn{1}{c|}{\rd(K)} & \multicolumn{1}{c|}{n} & \qquad &
\gms_2\text{ for }\KK&
\multicolumn{1}{c|}{\rd(K)} & \multicolumn{1}{c|}{n}\\
\cline{1-3} \cline{5-7}
4.002 & 16.032 & 88 &&
4.303 & 19.742 & 400 \\
4.066 & 16.756 & 110 &&
4.428 & 21.535 & 2400 \\
4.216 & 18.597 & 220 &&
4.449 & 21.843 & 4800 \\
4.231 & 18.788 & 240 &&
4.460 & 22.021 & 8862  \\
\cline{1-3} \cline{5-7}
\end{array}\]
\end{table}
These values are simply an extract from \cite{odlyzko-tables} and are
provided for easy reference.  We have added the values in the column
``$\gms_2\text{ for }\KK$'' which are simply log base $2$
of the values in the $\rd(K)$ column, and then rounded down.

\subsection{Degrees 9--11}
\label{deg9-11}
The main goal of this section is to prove Proposition~\ref{deg9-11thm}
below.  First, we establish some preliminaries.

\begin{prop} \label{gmsbnd11}
  If $K$ is an extension of $\Q$ of degree $n<12$, and $m$ is the
  number of wild slopes for $p=2$ for $K^\g\otimes\QQ$, then
  \[
  \gms_2(K) \leq \left\lbrace \begin{array}{@{}c@{\,}ll}
    97/24 & < 4.042 & \text{if $m\leq 4$} \\
    101/24 & < 4.209 & \text{if $m\leq 5$} \\
    53/12 &  < 4.417 & \text{if $m\leq 6$} \\
    71/16 & < 4.438 & \text{for any $m$} \\
  \end{array} \right.
\]
\end{prop}
\begin{proof}
  We consider the possible decompositions of $K\otimes\QQ\cong \prod_i
  K_{p,i}$.  If no $K_{p,i}$ has degree $8$ over $\QQ$, then all slopes $s$ for
  $K^\g\otimes\QQ$ satisfy $s\leq 4$ by
  \cite{jr-local-database}.  Hence, $\gms_2(K)\leq
  4$ which implies the asserted bounds.

  Now suppose some $K_{p,i}/\QQ$ is an octic extension.  There can be at
  most one other non-trivial extension among the $K_{p,i}/\QQ$, and
  its degree over $\QQ$ is at most $3$.  A complete summary of all
  candidates of the slope content of an
  octic over $\QQ$ is given in \cite{jr-2adicoctics}.
  Table~\ref{octicslopes} gives maximal slope content for $m$ wild
\begin{table}[htb]
\caption{Maximum \label{octicslopes}slope combinations for octic extensions of $\QQ$.}
\[
\begin{array}{|c|c|c|}
\hline
\text{\textbf{\# slopes}} & \text{\textbf{Slope Content}} & \gms_2 \\
\hline
3 & [3, 4, 5]_1 & 31/8 \\
4 & [2, 3, 4, 5]_1 & 4 \\
5 & [2, 3, 7/2, 4, 5]_1 & 67/16 \\
6 & [2, 3, 7/2, 4, 17/4, 5]_1 & 141/32 \\
\hline
\end{array}
\]
\end{table}
slopes, for $m\ge 3$.  Using this, we compute the crude bound for an
octic and a quadratic (maximal slope content being $[3]_1$) and for an
octic with a cubic (maximal slope content being $[\,]_3$).  The
statement of the theorem lists the resulting values of the Galois mean
slope.

For example, the first entry arises from the maximum contribution by
an octic with $3$ slopes, plus a 
single slope of $3$ from a quadratic to give slope content
$[3,3,4,5]_1$.  On the other hand, the maximum for $5$ slopes arises from an
octic with slope content $[2,3,7/2,4,5]_1$ and a tame cubic to give
content $[2,3,7/2,4,5]_3$.
\end{proof}
If $K/\Q$ is unramified away from $2$, then we can compare
$\gms_2(K)$ with values in Table~\ref{uncond} to get the following.
\begin{cor} \label{cor9-11}
  If $K/\Q$ is unramified away from $2$, $[K : \Q]<12$, and
  $G=\gal(K^\g/\Q)$, let $m=\ord(|G|)$, then 
\[ |G| < \begin{cases}
 110 & \text{if $m\leq 4$} \\
 220 & \text{if $m\leq 5$} \\
 2400 & \text{if $m\leq 6$} \\
 4800 & \text{in all cases} \\
\end{cases}
\]
\end{cor}

\begin{prop} \label{deg9-11thm}
  If $K\in\KK$, then $[K :  \Q]$ is not equal to $9$, $10$, or
  $11$.
\end{prop}
\begin{proof}
  We consider each of the possible Galois groups $G$ of polynomials of
  degree $9$, $10$, and $11$, of which there are $34$, $45$, and $8$
  groups respectively.  By Theorem~2.23 of \cite{harb}, we can
  eliminate $G$ if $|G|$ is not a multiple of $16$.  By
  Corollary~\ref{cor9-11}, we eliminate groups where $|G|\geq 4800$.
  Next, we eliminate groups which have a quotient which has already
  been eliminated.  
Note, this already eliminates all groups in degree $11$.
  Each of the remaining groups is then eliminated by
  Corollary~\ref{cor9-11} by comparing $|G|$ with $\ord(|G|)$:
\begin{align*}
|9T_{19}| &= 144=2^4 \cdot 9 & |10T_{28}| &= 400 = 2^4\cdot 25 &
|10T_{30}| &=720 = 2^4\cdot 45 \\ 
|10T_{31}| &= 720 = 2^4\cdot 45 &
|10T_{33}| &= 800 = 2^5\cdot 25 &
|10T_{35}| &= 1440 = 2^5\cdot 45 &
\end{align*}
\end{proof}

\subsection{Degrees 12--15}
\label{deg12}

The structure of this section is similar to that of
section~\ref{deg9-11}, although bounding $\gms_2$ is more
complicated.  We start with a bound on the slope content of composita
of certain quartic extensions. 

\begin{prop} \label{allquartic} %
  Let $F$ be the compositum of all quartic extensions of $\QQ$ whose
  Galois closures have Galois groups which are $2$-groups.  Then
  $[F:\QQ] = 2^8$, $F$ has residue field degree $4$
  and slope content $[2,2,3,3,7/2,4]_1$.
\end{prop}
\begin{proof}
  Clearly the tame degree is $1$ since the compositum has Galois group
  a $2$-group.  Let $\mathbf{G}_2$ be the Galois group of the
  composita of all $2$-group extensions of $\QQ$.  The group
  $\mathbf{G}_2$ is the pro-2 completion of the group with
  presentation $\langle x,\, y,\, z\mid x^2y^3z^{-1}yz=1\rangle$
  \cite{wingberg-book}.  Using this description, one can compute with
  \gap{} \cite{gap} the intersection of the kernels of all
  homomorphisms to the groups $V_4$, $C_4$, and $D_4$, the three
  Galois groups of quartics which are $2$-groups.  The quotient of
  $\mathbf{G}_2$ by this kernel has order $2^8$, hence $[F:\QQ]=2^8$.

  Naturally, this compositum contains the unramified extension of
  $\QQ$ of degree $4$, and 
  from the tables in \cite{jr-local-database}, we see that the wild
  slopes include $[2,2,3,7/2,4]$ since there are $D_4$ quartic fields
  with at least each slope once, and one with two slopes of $2$.

  To find the final slope, we consider the group $C_2:C_4:(C_4\times
  C_2)=\gggg$, meaning group number $61$ among groups of order $64$ in
  the numbering of \gap{}.  From the presentation above,
  one can check that $\gggg$ appears as a Galois group over $\QQ$.  
  From the group itself, one can verify that a field with Galois group
  $\gggg$ is the compositum of its $D_4$ subfields. Hence, there is an
  extension of $\QQ$ with Galois group $\gggg$ which is a subfield of
  $F$. But, the group
  $\gggg$ has $8T_{11}$ as a quotient.  Consulting
  \cite[Table~5.1]{jr-2adicoctics}, we see that there are $8T_{11}$
  fields with slope content $[2,3,3]_1$ in the notation used here (it is
  listed there as $[0,2,3,3]$).  In
  particular, there are two slopes of $3$ for $F$.
\end{proof}

\begin{rem}
  One can see the two slopes of $3$ explicitly as follows.  Consider
  the polynomials $x^4 + 2x^2 - 2$, $x^4 + 6x^2 + 3$, and $x^4 + 6x^2
  + 18$, which each have Galois group $D_4$ both over $\QQ$ and over
  $\Q$.  One can compute using \gp{} \cite{gp} their compositum
  over $\Q$, $K_{64}$, which is a degree $64$ extension with
  discriminant $2^{196}$.  The extension $K_{64}$ also has a single
  prime above $2$, so its global Galois group equals its decomposition
  group for the prime above $2$.  As a result, all subfields of the
  $2$-adic field are seen globally.  Computing subfields of $K_{64}$
  and the $2$-parts of their discriminants shows that $K_{64}$
  contains a quadratic unramified extension and has slope content
  $[2,2,3,3,7/2]_1$.
\end{rem}

\begin{prop} \label{deg15gms}
  If $K$ is an extension of $\Q$ of degree $n<16$, and $m$ is the
  number of wild slopes for $p=2$ for $K^\g\otimes\QQ$, then
  \[
  \gms_2(K) \leq \left\lbrace \begin{array}{@{}c@{\,}ll}
    203/48 & < 4.230 & \text{if $m\leq 4$} \\ 
    413/96 & < 4.303 & \text{if $m\leq 5$} \\ 
    495/112 & < 4.420 & \text{if $m\leq 6$} \\  
    107/24 & < 4.459 & \text{for any $m$} \\ 
  \end{array}\right.
\]
\end{prop}
\begin{proof}
  As in Proposition~\ref{gmsbnd11}, the local algebra $K\otimes \QQ =
  \prod_i K_{p,i}$ must have an octic field $K_{p,i}$ or all slopes would be
  $\leq 4$, here using
  Proposition~\ref{slopebound12} to rule out local fields $K_{p,i}$ with
  $9\leq [K_{p,i} : \QQ]\leq 15$.

  Note that a degree $6$ field can always be replaced with its twin
  algebra.  From \cite{web-local-database}, all $2$-adic sextic fields
  have twin algebras which split as a product of fields of degrees
  less than or equal to $4$.  Hence, we do not need to consider sextic
  factors.


  The cases with $m\leq 6$ work just like in
  Proposition~\ref{gmsbnd11}, where we use the crude bound for the
  slope content of the composita.  For example, our bound for $5$
  slopes comes from $[3,4,5]_1$ for the octic, $[3,4]_1$ for a
  quartic, and a tame cubic combining to yield
  $\gms_2([3,3,4,4,5]_3)=413/96$.

  For $m\geq 7$, we divide into several cases.  If $5$ is not a slope of
  the octic factor, we can apply the crude bound for the maximum slope
  content for the compositum of an octic (if $5$ is not a slope,
  $[3,7/2,4,17/4,19/4]_1$ has the largest $\gms_2$), a quartic with slope
  content $[2,3,4]_1$, and a quadratic
  with slope content $[3]_1$.  The result is $\gms_2([2$, $3$, $3$,
  $3$, $7/2$, $4$, $4$, $17/4$, $19/4]_1) = 421/96<107/24$.

  Now assume that $5$ is a slope for the octic.
  If $17/4$ is {\em not} a  slope of the octic factor, then the
  maximum slope content of the octic is $[2, 3, 7/2, 4, 5]_1$.  Again,
  the crude bound for this with a quartic and a quadratic is
  $\gms_2([2$, $2$, $3$, $3$, $3$, $7/2$, $4$, $4$, $5]) = 1125/256<107/24$. 

  Finally, we have the case where $5$ and $17/4$ are both slopes of
  the octic.  From \cite{jr-2adicoctics}, the slope content of such an
  octic is $[2,3,7/2,4,17/4,5]_1$ and only possibilities for
  the Galois group are $8T_{27}$, $8T_{28}$, are $8T_{35}$, each of
  which is a $2$-group.  In each
  case, the bottom $4$ slopes $[2,3,7/2,4]$ are visible in the
  compositum of quartic subfields.

  If we combine with quartics whose Galois groups are $2$-groups, then
  the maximal slope content of the quartic part is $[2,2,3,3,7/2,4]$
  by Proposition~\ref{allquartic},
  so maximum combination in this case is
  $\gms_2([2,2,3,3,7/2,4,17/4,5]_3)=427/96<107/24$.
  Finally, if we use the crude bound for the composita of an octic
  with slope content $[2,3,7/2,4,17/4,5]_1$, a
  quartic whose Galois group is not a 2-group, so maximal slope
  content of $[8/3,8/3]_3$, and a
  quadratic (slope content $[3]_1$), we get $\gms_2([2$, $8/3$, $8/3$,
  $3$, $3$, $7/2$, $4$, $17/4$, $5]_3)=107/24$.
\end{proof}

Now, we can combine Proposition~\ref{deg15gms} with bounds from
Table~\ref{uncond} to get the following.
\begin{cor} \label{cor12-15}
  If $K/\Q$ is unramified away from $2$, $[K : \Q]<16$, and
  $G=\gal(K^\g/\Q)$, let $m=\ord(|G|)$, then 
\[ |G| < \begin{cases}
 240 & \text{if $m\leq 4$} \\
 400 & \text{if $m\leq 5$} \\
 2400 & \text{if $m\leq 6$} \\
 8862 & \text{in all cases} \\
\end{cases}
\]
\end{cor}

\begin{prop} \label{deg12-15thm}
  If $K\in\KK$, then $[K : \Q]$ is not equal to $12$, $13$, $14$, or
  $15$.
\end{prop}
\begin{proof}
  As before, we consider each of the possible Galois groups $G$ of
  polynomials of the stated degrees.  For $n=12$, $13$, $14$, and
  $15$, there are $301$, $9$, $63$, and $104$ conjugacy classes of
  subgroups in $S_n$ respectively.  By Theorem~2.23 of \cite{harb},
  we can eliminate $G$ if $|G|$ is not a multiple of $16$.  By
  Corollary~\ref{cor12-15}, we eliminate groups where $|G|\geq 8862$,
  and then eliminate groups which have a quotient which has already
  been eliminated.  For the remaining groups, we give their orders
  with partial factorization to show that they too are eliminated by
  \ref{cor12-15}.
\[ \setlength{\arraycolsep}{1pt}
\begin{array}{rrlcl}
|12T_{j}| =& 1296=& 2^4 \cdot 81 &\qquad\qquad& \text{for $j=215,\ 216$} \\
|12T_{j}| =& 2592=& 2^5 \cdot 81 && \text{for $244\leq j\leq 249$} \\
|12T_{j}| =& 5184=& 2^6 \cdot 81 && \text{for $262\leq j\leq 264$} \\
|13T_{7}| =& 5616 =& 2^4\cdot 351 && \\
|14T_{16}| =& 336 =& 2^4\cdot 21 &&
\end{array}
\]
%
\end{proof}

Note, no transitive subgroups of $S_{15}$ passed through the various
filters discussed in the proof of Theorem~\ref{deg12-15thm}, and only
one group needed to be considered in each of degrees $13$ and $14$.


Since there is particular interest in whether or not simple groups are
in $\G$, we extract the new cases covered by Theorem~\ref{mainthm}.
Additional results on simple groups excluded from $\G$ by a
combination of root discriminant bounds and group theoretic
techniques, see \cite{jj-simple}.

\begin{cor}\label{simple2cor}
  The following simple groups are not elements of $\G$: alternating groups
  $A_j$ for $9\leq j \leq 15$, Mathieu groups $M_{11}$ and
  $M_{12}$, and $\PSL_3(3)$.
\end{cor}

\bibliographystyle{alpha}
\bibliography{awayfrom2}

\end{document}